\documentclass{amsart}
\usepackage{multicol,graphicx,color}
\usepackage{amsthm}
\usepackage{amsmath}
\usepackage{amssymb}
\usepackage{amsfonts} 

\theoremstyle{plain}
\newtheorem{theorem}{Theorem}[section]

\newtheorem{lemma}[theorem]{Lemma}

\begin{document}

\title[Periodic SBC Orbits in the PPS4BP]{Existence and Stability of Symmetric Periodic Simultaneous Binary Collision Orbits in the Planar Pairwise Symmetric Four-Body Problem}
\author[Bakker]{Lennard F. Bakker}
\author[Ouyang]{Tiancheng Ouyang} 
\author[Yan]{Duokui Yan}
\author[Simmons]{Skyler Simmons}
\address{Department of Mathematics \\  Brigham Young University\\ Provo, UT USA 84602}
\email[Lennard F. Bakker]{bakker@math.byu.edu}
\email[Tiancheng Ouyang]{ouyang@math.byu.edu}
\email[Skyler Simmons]{xinkaisen@gmail.com}
\address{Chern Institute of Mathematics \\ Nankai University \\ Tianjin 300071, P.R.China}
\email[Duokui Yan]{duokuiyan@gmail.com}


\begin{abstract} We extend our previous analytic existence of a symmetric periodic simultaneous binary collision orbit in a regularized fully symmetric equal mass four-body problem to the analytic existence of a symmetric periodic simultaneous binary collision orbit in a regularized planar pairwise symmetric equal mass four-body problem. We then use a continuation method to numerically find symmetric periodic simultaneous binary collision orbits in a regularized planar pairwise symmetric 1, m, 1, m four-body problem for $m$ between $0$ and $1$. Numerical estimates of the the characteristic multipliers show that these periodic orbits are linearly stability when $0.54\leq m\leq 1$, and are linearly unstable when $0<m\leq0.53$.
\end{abstract}

\maketitle

\section{Introduction}

In the $N$-body problem, linearly stable periodic orbits may trap around themselves bounded, non-chaotic motion of the $N$ masses \cite{Co}. Some of the known examples of linearly stable periodic orbits in the three-body problem are the elliptic Lagrangian triangular periodic orbits for certain values of eccentricity and the three masses  \cite{MeSc},\cite{Ro1}, and the Montgomery-Chenciner figure-eight periodic orbit  for equal masses \cite{CM},\cite{MR},\cite{Ro2},\cite{HS1},\cite{HS2}. Other examples of linearly stable periodic orbits in the three or four-body problem involve binary collisions (BC) and/or simultaneous binary collisions (SBCs). The regularization of these kinds of singularities is achieved by a generalized Levi-Civita type transformation and an appropriate scaling of time, as adapted from Aarseth and Zare \cite{AZ} to the particular problem (also see \cite{MS},\cite{OX}).

Schubart \cite{Sc} numerically discovered a singular symmetric periodic orbit in the collinear three-body equal mass problem. In this orbit, the inner mass alternates between binary collisions with the two outer masses. H\'enon \cite{He} extended Schubart's numerical investigations to the case of unequal masses. Only recently did Venturelli \cite{Ve} and Moeckel \cite{Mo} prove the analytic existence of the Schubart orbit when the outer masses are equal and the inner mass is arbitrary. The linear stability of the Schubart orbit was determined numerically by Hietarinta and Mikkola \cite{HM} revealing that linear stability occurs only for certain choices of the three masses. Numerically, non-Schubart-like linearly stable periodic orbits in the collinear three-body problem were found by Saito and Tanikawa for certain choices of the masses \cite{ST1},\cite{ST2},\cite{ST3}.

Sweatman \cite{SW}, \cite{SW2} and Sekiguchi and Tanikawa \cite{SeT} numerically found a symmetric  Schubart-like orbit in the symmetric collinear four-body problem with masses $1$, $m$, $m$, and $1$. This Schubart-like periodic orbit alternates between SBCs of the two outer pairs of masses and binary collisions of the inner two masses. Ouyang and Yan \cite{OY2} proved analytically the existence and symmetry of this orbit. In the regularized setting, this periodic orbit has a symmetry group isomorphic to $D_2$, of which both of the generators are time-reversing symmetries. (Here the dihedral group $D_k$ is the group of symmetries of the regular $k$-gon.) Sweatman \cite{SW2} numerically showed that the Schubart-like orbit is linearly stable when $0<m<2.83$ or $m>35.4$, and is otherwise unstable. This linear stability was confirmed \cite{BOYS} using Robert's symmetry reduction technique \cite{Ro2}.

Ouyang, Yan, and Simmons \cite{OY3} numerically found and analytically proved the existence and symmetries of a singular symmetric periodic orbit in the fully symmetric planar four-body problem with equal masses. (In the fully symmetric planar four-body equal mass problem, the position of one mass determines the positions of the other three masses.) In this orbit, the four masses alternate between different SBC's. In the regularized setting, this periodic orbit has a symmetry group isomorphic to $D_4$, of which one of the two generators is a time-reversing symmetry, while the other generator is a time-preserving symmetry. By using Roberts' symmetry reduction method, we showed \cite{BOYS} that this symmetric periodic simultaneous binary collision orbit is linearly stable.

In this paper we extend the analytic existence of a symmetric periodic SBC orbit in the fully symmetric planar four-body equal mass problem \cite{OY3} to the analytic existence of a symmetric periodic SBC orbit in the planar pairwise symmetric four-body problem, or PPS4BP for short. The positions of the four pairwise symmetric bodies in the plane are $(x_1,x_2)$, $(x_3,x_4)$, $(-x_1,-x_2)$, and $(-x_3,-x_4)$, where the corresponding masses are $1$, $m$, $1$, $m$ with $0<m\leq 1$. With $t$ as the time variable and $\dot{}=d/dt$, the momenta for the four masses are $(\omega_1,\omega_2) = 2(\dot x_1,\dot x_2)$, $(\omega_3,\omega_4)=2m(\dot x_3,\dot x_4)$, $-(\omega_1,\omega_2)$, and $-(\omega_3,\omega_4)$. The Hamiltonian for the PPS4BP is $H=K-U$, where
\[ K = \frac{1}{4}\big[ \omega_1^2 + \omega_2^2\big] + \frac{1}{4m}\big[ \omega_3^2 + \omega_4^2\big],\]
and
\begin{align*} U = & \frac{1}{2\sqrt{x_1^2 + x_2^2}} + \frac{2m}{\sqrt{(x_3-x_1)^2 + ( x_4-x_2)^2}}\\
& + \frac{2m}{\sqrt{(x_1+x_3)^2+(x_2+x_4)^2 }} + \frac{m^2}{2\sqrt{x_3^2 + x_4^2}}.
\end{align*}
The angular momentum for the PPS4BP is
\[ A = x_1\omega_2-x_2\omega_1 + x_3\omega_4- x_4\omega_3.\]
The center of mass is fixed at the origin, and the linear momentum is zero. With
\[ J = \begin{bmatrix} 0 & I \\ -I & 0\end{bmatrix}\]
for $I$ the $4\times 4$ identity matrix, the vector field for the PPS4BP is $J\nabla H$, i.e., the Hamiltonian system of equations with Hamiltonian $H$ are $\dot x_i = \partial H/\partial \omega_i$, $\dot \omega_i = -\partial H/\partial x_i$, $i=1,2,3,4$. The PPS4BP presented here is the Caledonian symmetric four-body problem \cite{SSS} with non-collinear initial positions.

\begin{figure}\scalebox{0.35}{\includegraphics{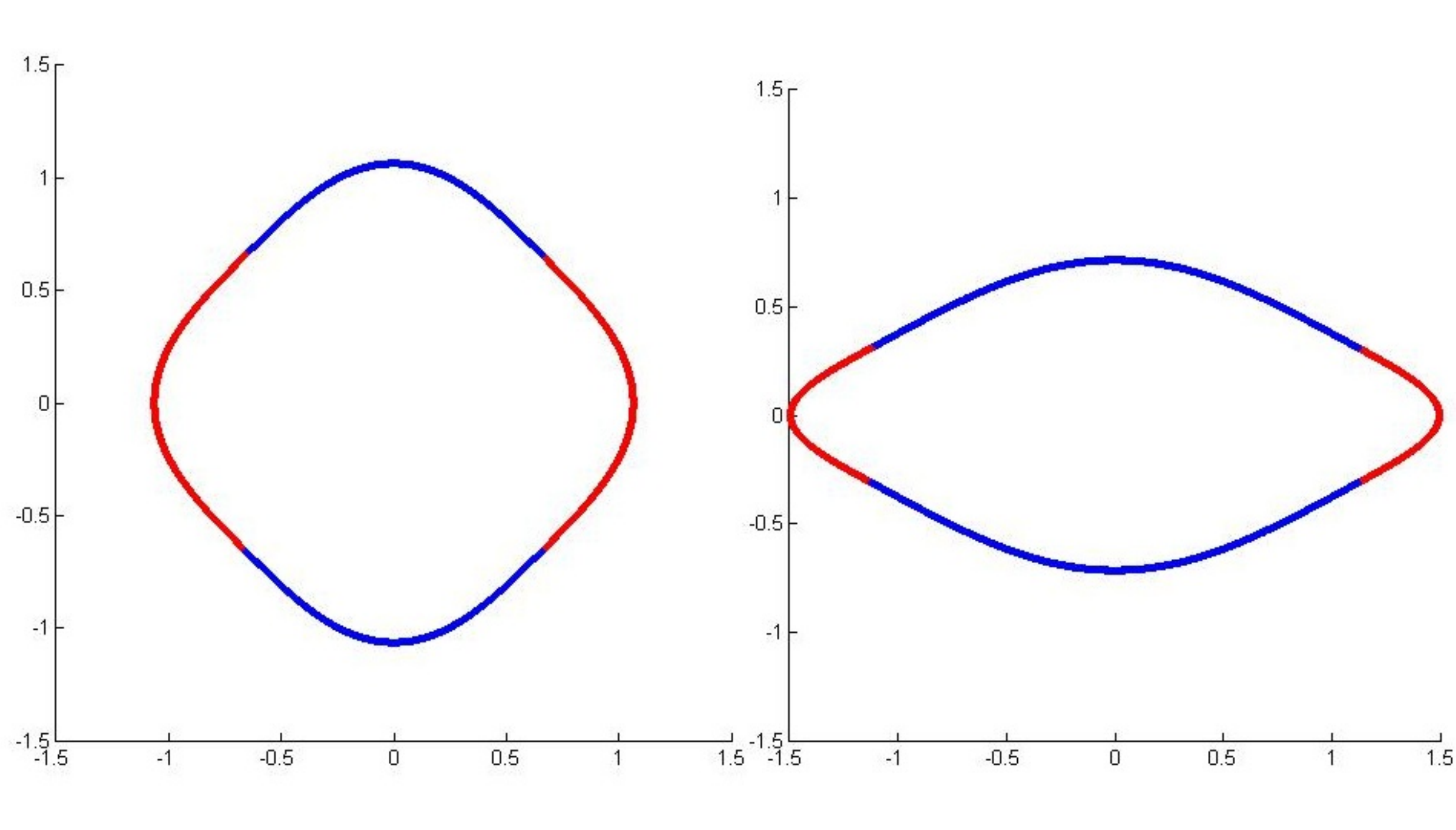}}\caption{The symmetric periodic SBC orbit in the PPS4BP for $m=1$ (left) and $m=0.539$ (right). The two red curves are those traced out by $\pm(x_1(t),x_2(t))$, and the two blue curves are those traced out by $\pm(x_3(t),x_4(t))$.}\label{figure1}\end{figure}

The initial conditions for the orbits of interest has the first body of mass $1$ located on the positive  horizontal axis with its momentum perpendicular to the horizontal axis, and the first body of mass $m$ located on the positive vertical axis with its momentum perpendicular to the vertical axis. Specifically, at $t=0$ we have
\begin{align*} &x_1 >0,\ x_2=0,\ x_3=0,\ x_4>0, {\rm\ with\ } x_4\leq x_1, \\
&\omega_1=0,\ \omega_2>0,\ \omega_3>0,\ \omega_4=0, {\rm\ with\ } \omega_2\leq\omega_3,
\end{align*}
at which $H$ is defined. The first objective is to find, for $0<m\leq 1$, values of $x_1,x_4,\omega_2,\omega_3$ at $t=0$ such that (i) $x_3-x_1=0$ and $x_4-x_2=0$ with $x_1^2+x_2^2\ne0 $ at some $t=t_0>0$, (ii) $x_1+x_3=0$ and $x_2+x_4=0$ with $x_1^2+x_2^2\ne0$ at some $t=t_1>t_0$, (iii) the orbit extends to a symmetric periodic orbit, and (iv) the periodic orbit avoids all the other kinds of collisions. Such an orbit experiences a SBC in the first and third quadrant at $t=t_0$, and then another SBC in the second and fourth quadrants at $t=t_1$, before returning to its initial conditions at some $t=t_2>t_1$. The presence of collisions along the orbit necessarily imposes zero angular momentum on the orbit, thus requiring that $x_1\omega_2-x_4\omega_3=0$ at $t=0$. Examples of these symmetric periodic SBC orbits in the PPS4BP with masses $1$, $m$, $1$, $m$ are illustrated in Figure \ref{figure1} for $m=1$ and $m=0.539$. The second objective is to numerically investigate the linear stability of the symmetric periodic SBC orbits as $m$ varies over interval $(0,1]$.

The regularization of the SBCs, as described by (i) and (ii) above, in the Hamiltonian system of equations with Hamiltonian $H$ plays a key role in achieving the two objectives. Section 2 details this regularization which consists of two canonical transformations followed by a scaling of time $t=\theta(s)$ with $s$ as the regularizing time variable, producing a new Hamiltonian $\hat\Gamma$ for the PPS4BP in extended phase space. Section 3 describes a scaling of orbits of the Hamiltonian system of equations with Hamiltonian $\hat\Gamma$ which shows that any such periodic solution always belongs to a one-parameter family of periodic solutions for which the linear stability is the same for all periodic solutions in the family. Section 4 describes the symmetries of the Hamiltonian $\hat\Gamma$ which are used to construct periodic solutions with a $D_4$ symmetry group generated by a time-reversing symmetry and a time-preserving symmetry.

In Sections 5 and 6, we prove the analytic existence of a periodic SBC orbit $\gamma(s)$, with a $D_4$ symmetry group, for the Hamiltonian system of equations with Hamiltonian $\hat\Gamma$ with $m=1$, and numerically investigate its stability. The proof extends the analytic existence of a symmetric periodic SBC orbit in the fully symmetric planar four-body equal mass problem, as found in \cite{OY3}, to the PPS4BP with equal masses. Our numerical estimates of the characteristic multipliers show that this periodic orbit is linearly stable.

In Section 7, we numerically continue the linearly stable symmetric periodic SBC orbit $\gamma(s)$ from $m=1$ to a symmetric periodic SBC orbit $\gamma(s;m)$ for $m<1$, and then investigate the linear stability of $\gamma(s;m)$ as $m$ varies in the interval $(0,1]$. We use trigonometric polynomials as approximations of the periodic orbits (cf.\,\cite{Si}). The numerical algorithm for continuation starts with a trigonometric polynomial approximation of $\gamma(s;1)=\gamma(s)$ that is used as a initial guess for a trigonometric polynomial approximation of $\gamma(s;0.99)$, where the coefficients of the trigonometric polynomial are optimized through a variational approach. This process is repeated, using the optimized approximation of $\gamma(s;0.99)$ as the initial guess for $\gamma(s;0.98)$, etc., until an optimized approximation of $\gamma(s;0.01)$ is obtained. Numerical estimates of the characteristic multipliers of $\gamma(s;m)$ show it is linearly stable when $0.54\leq m\leq 1$ and  unstable when $0<m\leq 0.53$. We have analyzed the linear stability of these orbits for $0<m\leq 1$, especially for $0.53<m<0.54$ with increments of $0.001$, using Roberts' symmetry reduction method, the results of which will appear in a subsequent paper \cite{BOSM}.

\section{Regularization} We adapt the regularization of Aarseth and Zare \cite{AZ} to the PPS4BP to regularize SBCs as described in the first objective. This regularization differs from the one used in the Caledonian symmetric four-body problem \cite{SSS} in that we only regularize the SBCs as described in the first objective of the Introduction.

The first canonical transformation in our regularization is
\[ (x_1,x_2,x_3,x_4,\omega_1,\omega_2,\omega_3,\omega_4) \to (g_1,g_2,g_3,g_4,h_1,h_2,h_3,h_4)\]
determined by the generating function
\[ F_1(x_1,x_2,x_3,x_4,h_1,h_2,h_3,h_4) = h_1(x_1-x_3) + h_2(x_2-x_4) + h_3(x_1+x_3) + h_4(x_2+x_4).\]
So the first canonical transformation is determined by
\begin{equation}\label{firstcantrans} \omega_i = \frac{\partial F_1}{\partial x_i}, \ \ g_i = \frac{\partial F_1}{\partial \omega_i}, \ \ i=1,2,3,4.\end{equation}
The new Hamiltonian is $\tilde H=\tilde K - \tilde U$, where
\[ \tilde K = \frac{ (h_1+h_3)^2+(h_2+h_4)^2}{4} + \frac{(h_3-h_1)^2 + (h_4-h_2)^2}{4m},\]
and
\[ \tilde U = \frac{1}{\sqrt{(g_1+g_3)^2 + (g_2+g_4)^2}} + \frac{2m}{\sqrt{g_1^2+g_2^2}} + \frac{2m}{\sqrt{g_3^2+ g_4^2}} + \frac{m^2}{\sqrt{(g_1-g_3)^2 + (g_2-g_4)^2}}.\]

The second canonical transformation in our regularization,
\[ (g_1,g_2,g_3,g_4,h_1,h_2,h_3,h_4) \to (u_1,u_2,u_3,u_4,v_1,v_2,v_3,v_4)\]
is determined by the generating function
\[ F_2 (h_1,h_2,h_3,h_4,u_1,u_2,u_3,u_4) = - \sum_{j=1}^4 h_j f_j(u_1,u_2,u_3,u_4),\]
where
\[ f_1 = u_1^2-u_2^2,\ \ f_2 = 2u_1u_2, \ \ f_3 = u_3^2 - u_4^2, \ \ f_4 = 2u_3u_4.\]
So the second canonical transformation is determined by
\begin{equation}\label{secondcantrans} g_i = -\frac{\partial F_2}{\partial h_i}, \ \ v_i = - \frac{\partial F_2}{\partial u_i}, \ \ i=1,2,3,4.\end{equation}
The new Hamiltonian is $\hat H = \hat K - \hat U$, where
\begin{align*}
\hat K = & \frac{1}{16}\bigg(1+\frac{1}{m}\bigg)\bigg[ \frac{(v_1^2+v_2^2)(u_3^2+u_4^2)  + (v_3^2+v_4^2)(u_1^2+u_2^2)}{(u_1^2+u_2^2)(u_3^2+u_4^2)} \bigg] \\
& + \frac{1}{8}\bigg(1-\frac{1}{m}\bigg) \frac{(v_3u_3-v_4u_4)(v_1u_1-v_2u_2) + (v_3u_4+v_4u_3)(v_1u_2+v_2u_1)}{(u_1^2+u_2^2)(u_3^2+u_4^2)},
\end{align*}
and
\begin{align*} \hat U = & \frac{1}{\sqrt{(u_1^2-u_2^2 +u_3^2 - u_4^2)^2 + (2u_1u_2+2u_3u_4)^2}} + \frac{2m}{u_1^2+u_2^2} + \frac{2m}{u_3^2+u_4^2} \\
& + \frac{m^2}{\sqrt{(u_1^2-u_2^2-u_3^2+u_4^2)^2 + (2u_1u_2-2u_3u_4)^2 }}.
\end{align*}

We introduce a new time variable $s$  by the regularizing change of time
\[ \frac{dt}{ds} = (u_1^2+u_2^2)(u_3^2+u_4^2).\]
To simplify notation, we set
\begin{align*}
M_1 & = v_1u_1-v_2u_2, & M_2 & = v_1u_2+v_2u_1, \\
M_3 & = v_3u_3-v_4u_4, & M_4 & = v_3u_4 + v_4u_3, \\
M_5 & = u_1^2 -u_2^2 + u_3^2 - u_4 ^2, & M_6 & = 2u_1u_2+2u_3u_4, \\
M_7 & = u_1^2-u_2^2-u_3^2 + u_4^2, & M_8 & = 2u_1u_2-2u_3u_4.
\end{align*}
The Hamiltonian in the extended phase space with coordinates $u_1$, $u_2$, $u_3$, $u_4$, $\hat E$, $v_1$, $v_2$, $v_3$, $v_4$, $t$ is
\begin{align*} \hat\Gamma =  \frac{dt}{ds}\big(\hat H - \hat E\big) 
= & \frac{1}{16}\bigg(1+\frac{1}{m}\bigg)\bigg( (v_1^2+v_2^2)(u_3^2+u_4^2) + (v_3^2+v_4^2)(u_1^2+u_2^2) \bigg) \\
& +  \frac{1}{8}\bigg(1-\frac{1}{m}\bigg) \big(M_3M_1 + M_4M_2\big) \\
& - \frac{(u_1^2+u_2^2)(u_3^2+u_4^2)}{\sqrt{M_5^2 + M_6^2}} - 2m\big(u_1^2+u_2^2+u_3^2+u_4^2) \\
& - \frac{m^2(u_1^2+u_2^2)(u_3^2+u_4^2)}{\sqrt{M_7^2 + M_8^2 }} - \hat E(u_1^2+u_2^2)(u_3^2+u_4^2).
\end{align*}
With ${}^\prime=d/ds$, the Hamiltonian system of equations with Hamiltonian $\hat\Gamma$ is
\begin{equation}\label{regularized} u_i^\prime  = \frac{\partial \hat\Gamma}{\partial v_i}, \ \ v_i^\prime = -\frac{\partial \hat\Gamma}{\partial u_i}, \ \ i=1,2,3,4,\end{equation}
along with the auxiliary equations,
\begin{equation}\label{timechange} \hat E^\prime = \frac{\partial\hat\Gamma}{\partial t} =  0,\ \ t^\prime  = -\frac{\partial\hat\Gamma}{\partial \hat E} = (u_1^2+u_2^2)(u_3^2+u_4^2).\end{equation}

On the level set $\hat\Gamma=0$, the value of the Hamiltonian $\hat H$ (i.e., the energy) along solutions of the Hamiltonian system of equations with Hamiltonian $\hat \Gamma$ is $\hat E$. Independent of the values of $\hat  E$ and $\hat\Gamma$, the angular momentum $A=x_1\omega_2-x_2\omega_1+x_3\omega_4-x_4\omega_3$ in the coordinates $u_1,u_2,u_3,u_4,v_1,v_2,v_3,v_4$ simplifies to
\begin{equation}\label{ram} A = \frac{1}{2}\big[ -v_1u_2 + v_2u_1 - v_3u_4 + v_4u_3\big].\end{equation}

On the level set $\hat\Gamma=0$, two simultaneous binary collisions in the PPS4BP have been regularized in the Hamiltonian system of equations with Hamiltonian $\hat\Gamma$. The simultaneous binary collision $x_3-x_1=0$ and $x_4-x_2=0$ with $x_1^2+x_2^2\ne 0$ corresponds to $u_1^2+u_2^2=0$ with $u_3^2+u_4^2\ne 0$. These imply that $M_5^2+M_6^2=4(x_1^2+x_2^2)\ne0$ and $M_7^2+M_8^2=4(x_1^2+x_2^2)\ne0$. From $\hat\Gamma=0$ it follows that
\begin{equation}\label{rmomentum} v_1^2+v_2^2 = (32m^2)/(m+1).\end{equation}
Similarly, the simultaneous binary collision $x_3+x_1=0$ and $x_4+x_2=0$ with $x_1^2+x_2^2\ne 0$ corresponds to $u_3^2+u_4^2=0$ with $u_1^2+u_2^2\ne 0$, and hence that $M_5^2+M_6^2 = 4(x_1^2+x_2^2)\ne 0$, $M_7^2+M_8^2=4(x_1^2+x_2^2)\ne 0$, and, from $\hat\Gamma=0$, that
\[ v_3^2+v_4^2 = (32m^2)/(m+1).\]

On the level set $\hat\Gamma=0$, the other singularities of the PPS4BP have not been regularized in the Hamiltonian system of equations with Hamiltonian $\hat\Gamma$. The binary collision $x_1= 0$, $x_2=0$ with $x_3^2+x_4^2\ne0$ corresponds to $M_5^2+M_6^2=0$ and $M_7^2+M_8^2\ne 0$ with $u_1^2+u_2^2\ne 0$ and $u_3^2+u_4^2\ne 0$. The binary collision $x_3=0$, $x_4=0$ with $x_1^2+x_2^2\ne0$ corresponds to $M_5^2+M_6^2\ne0$ and $M_7^2+M_8^2=0$ with  $u_1^2+u_2^2\ne 0$ and $u_3^2+u_4^2\ne 0$. Because of the pairwise symmetry, there are no triple collisions. Total collapse $x_1=0$, $x_2=0$, $x_3=0$, $x_4=0$ corresponds to $u_1=0$, $u_2=0$, $u_3=0$, and $u_4=0$. A solution of the Hamiltonian system of equations with Hamiltonian $\hat\Gamma$ is called {\it nonsingular}\, if it avoids the unregularized binary collisions and total collapse singularities.

We establish next the correspondence between the original coordinates and the regularized coordinates for the initial conditions given in the Introduction, the proof of which is straight-forward.

\begin{lemma}\label{initial} The conditions $($at $t=0)$
\begin{align*} &x_1 >0,\ x_2=0,\ x_3=0,\ x_4>0, {\rm\ with\ } x_4\leq x_1, \\
&\omega_1=0,\ \omega_2>0,\ \omega_3>0,\ \omega_4=0, {\rm\ with\ } \omega_2\leq\omega_3,
\end{align*}
correspond to the conditions $($at $s=0)$
\begin{align*}
&u_3 = \pm u_1,\ u_4 =\mp u_2,{\rm\ with\ }u_1u_2<0,\ \vert u_2\vert\leq (\sqrt 2 - 1)\vert u_1\vert,\\
& v_3=\mp v_1,\ v_4 = \pm v_2,{\rm \ with\ } 0<v_1u_2+v_2u_1\leq v_2u_2-v_1u_1.
\end{align*}
\end{lemma}

\section{A Scaling of Periodic Orbits and Linear Stability}

A certain scaling of solutions of the Hamiltonian system of equations with Hamiltonian $\hat\Gamma$ produces more solutions. When applied to a periodic solution, this scaling leads to a one-parameter family of periodic solutions. The proof of the following result is a straight-forward verification.

\begin{lemma}\label{energy} If $\gamma(s) = (u_1(s),u_2(s),u_3(s),u_4(s),v_1(s),v_2(s),v_3(s),v_4(s))$ is a periodic solution of the Hamiltonian system of equations with Hamiltonian $\hat\Gamma$ on the level set $\hat\Gamma=0$ with period $T$ and energy $\hat E$, then for every $\epsilon>0$, the function
\[ \gamma_\epsilon(s) = (\epsilon u_1(\epsilon s),\epsilon u_2(\epsilon s),\epsilon u_3(\epsilon s),\epsilon u_4(\epsilon s),v_1(\epsilon s),v_2(\epsilon s),v_3(\epsilon s),v_4(\epsilon s))\]
is a periodic solution of the Hamiltonian system of equations with Hamiltonian $\hat\Gamma$ on the level set $\hat\Gamma=0$ with period $T_\epsilon=\epsilon^{-1}T$ and energy $\hat E_\epsilon=\epsilon^{-2} \hat E$.
\end{lemma}

The linear stability of a periodic orbit of the Hamiltonian system of equations with Hamiltonian $\hat\Gamma$ is determined by the linearization of the equations along the periodic orbit. By Lemma \ref{energy}, a periodic orbit $\gamma(s)$ on the level set $\hat\Gamma=0$ with period $T$ and energy $\hat E$ embeds into a one-parameter family $\gamma_\epsilon(s)$ of periodic orbits on the level set $\hat\Gamma=0$ with period $T_\epsilon$ and energy $\hat E_\epsilon$. The linearization of the Hamiltonian system of equations (\ref{regularized}) with Hamiltonian $\hat\Gamma$ along the periodic orbit $\gamma_\epsilon(s)$ is
\[ X^\prime = J\nabla^2\hat\Gamma(\gamma_\epsilon(s))X\]
where $\nabla^2\hat\Gamma$ is the matrix of second-order partials of $\hat\Gamma$. Let $X_\epsilon(s)$ be the solution of the linearization of the equations along $\gamma_\epsilon(s)$ that satisfies $X_\epsilon(0)=I$ (the $8\times 8$ identity matrix). The monodromy matrix for $\gamma_\epsilon(s)$ is $X_\epsilon(T_\epsilon)$, and the eigenvalues of $X_\epsilon(T_\epsilon)$ are the characteristic multipliers of $\gamma_\epsilon(s)$. A characteristic multiplier $\lambda$ of $\gamma_\epsilon(s)$ is defective if its geometric multiplicity is smaller than its algebraic multiplicity, i.e., its generalized eigenspace $\cup_{j\geq 1}{\rm ker}(X_\epsilon(T_\epsilon)-\lambda I)^j$ is not the same as its eigenspace ${\rm ker}(X_\epsilon(T_\epsilon)-\lambda I)$. Proving the following result is routine (see \cite{MH}).

\begin{lemma}\label{eigenvalueone} If $\gamma(s)$ is a periodic orbit of the Hamiltonian system of equations with Hamiltonian $\hat\Gamma$ on the level set $\hat\Gamma=0$, then for each $\epsilon>0$, the periodic orbit $\gamma_\epsilon(s)$ has $1$ as a defective characteristic multiplier with algebraic multiplicity at least two.
\end{lemma}

For each $\epsilon>0$, the periodic orbit $\gamma_\epsilon(s)$ is {\it spectrally stable}\, if all of its characteristic multipliers have modulus one.  By Lemma \ref{eigenvalueone}, the periodic orbit $\gamma_\epsilon(s)$ has $1$ as a defective characteristic multiplier with algebraic multiplicity at least two, and so the monodromy matrix $X_\epsilon(T_\epsilon)$ it not semisimple. However,  the two-dimensional subspace
\[ U_1 = {\rm Span}\left ( \gamma_\epsilon^{\, \prime}(0), \bigg(\frac{\partial}{\partial\epsilon}\gamma_\epsilon\bigg)(0)\right)\]
is $X_\epsilon(T_\epsilon)$-invariant. The periodic orbit $\gamma_\epsilon(s)$ is said to be {\it linearly stable}\, if it is spectrally stable and there exists a $6$-dimensional $X_\epsilon(T_\epsilon)$-invariant subspace $U_2$ such that $U_1+U_2={\mathbb R}^8$ and $X_\epsilon(T_\epsilon)$ restricted to $U_2$ is semisimple. A proof of the following is routine.

\begin{theorem}\label{stability} Suppose $\gamma(s)$ is a periodic orbit of the Hamiltonian system of equations with Hamiltonian $\hat\Gamma$ on the level set $\hat\Gamma=0$. Then $\gamma_\epsilon(s)$ is spectrally $($linearly$)$ stable for some $\epsilon>0$ if and only if $\gamma_\epsilon(s)$ is spectrally $($linearly$)$ stable for all $\epsilon>0$.
\end{theorem}

\section{Symmetries} The Hamiltonian system of equations with Hamiltonian $\hat\Gamma$ has a group of symmetries isomorphic to the dihedral group $D_4 = \langle a,b: a^2=b^4=(ab)^2=e\rangle$. With
\[ F = \begin{bmatrix} -1 & 0 \\ 0 & 1\end{bmatrix}, \ \ G= \begin{bmatrix}1 & 0 \\ 0 & 1\end{bmatrix},\]
define the matrices
\[ S_F = \begin{bmatrix} 0 & F & 0 & 0 \\ -F & 0 & 0 & 0 \\ 0 & 0 & 0 & F \\ 0 & 0 & -F & 0\end{bmatrix}, \ \ S_G = \begin{bmatrix} -G & 0 & 0 & 0 \\ 0 & G  & 0 & 0 \\ 0 & 0 & G & 0 \\ 0 & 0 & 0 & -G\end{bmatrix}.\]
These matrices satisfy $S_F^2=-I$, $S_F^4=I$, $S_G^2=I$, and $(S_FS_G)^2=I$.
Fixing the value of $\hat E$, these matrices satisfy $\hat\Gamma\circ S_F  = \hat \Gamma$ and $\hat\Gamma\circ S_G=\hat \Gamma$, and so $S_F$ and $S_G$ are the generators of the $D_4$-symmetry group for $\hat \Gamma$. If
\[ \gamma(s) = (u_1(s),u_2(s),u_3(s),u_4(s),v_1(s),v_2(s),v_3(s),v_4(s))\]
is a solution of the Hamiltonian system of equations with Hamiltonian $\hat\Gamma$, then $S_F\gamma(s)$, $S_F^2\gamma(s)$, and $S_G\gamma(-s)$ are also solutions of the Hamiltonian system of  equations with Hamiltonian $\hat\Gamma$. This means that $S_F$ is a time-preserving symmetry and that $S_G$ is a time-reversing symmetry. The proof of the following is routine.

\begin{lemma}\label{symmetry} If for some $s_0>0$ there is a nonsingular solution $\gamma(s)$, $s\in[0,s_0]$, of the Hamiltonian system of equations with Hamiltonian $\hat \Gamma$ such that for constants $\zeta_1\ne 0$, $\zeta_2\ne 0$, $\rho_1\ne 0$, and $\rho_2\ne 0$ there holds
\begin{align*} & u_1(0) = \zeta_1,\  u_2(0) = \zeta_2,\ u_3(0) = \zeta_1,\ u_4(0) = -\zeta_2, \\
& v_1(0) = \rho_1,\  v_2(0) = \rho_2, \ v_3(0) = -\rho_1, \ v_4(0) = \rho_2, 
\end{align*}
and
\begin{align*}
& u_1(s_0) = 0,\ u_2(s_0) = 0,\ u_3(s_0) \ne 0,\ u_4(s_0) \ne 0, \\
& v_1(s_0) \ne 0,\ v_2(s_0) \ne 0,\ v_3(s_0) =0,\ v_4(s_0) = 0, 
\end{align*}
then $\gamma(s)$ extends to a periodic orbit with period $8s_0$ and a symmetry group isomorphic to $D_4$ such that
\begin{align*}
& u_1(3s_0) \ne 0,\ u_2(3s_0) \ne 0,\ u_3(3s_0) = 0,\ u_4(3s_0) = 0, \\
& v_1(3s_0) = 0,\ v_2(3s_0) = 0,\ v_3(3s_0) \ne0,\ v_4(3s_0) \ne 0, 
\end{align*}
and
\begin{align*}
& u_1(5s_0) = 0,\ u_2(5s_0) = 0,\ u_3(5s_0) \ne 0,\ u_4(5s_0) \ne 0, \\
& v_1(5s_0) \ne 0,\ v_2(5s_0) \ne 0,\ v_3(5s_0) =0,\ v_4(5s_0) = 0, 
\end{align*}
and
\begin{align*}
& u_1(7s_0) \ne 0,\ u_2(7s_0) \ne 0,\ u_3(7s_0) = 0,\ u_4(7s_0) = 0, \\
& v_1(7s_0) = 0,\ v_2(7s_0) = 0,\ v_3(7s_0) \ne0,\ v_4(7s_0) \ne 0. 
\end{align*}
\end{lemma}

\section{Analytic Existence in the Equal Mass Case}
When $m=1$, there is an additional symmetry in the positions of the four masses that reduces the PPS4BP with equal masses to the fully symmetric planar four-body equal mass problem. We exploit this reduction to prove the existence of a symmetric periodic simultaneous binary collision orbit in the equal mass case.

The additional symmetry is the Ans\"atz, $x_4=x_1$, $x_3=x_2$, with $\vert x_2\vert \leq x_1$. From this it follows that $\omega_4 = \omega_1$, $\omega_3=\omega_2$, $x_1-x_2\geq 0$, $x_1+x_2\geq 0$. From the canonical transformations (\ref{firstcantrans}) and (\ref{secondcantrans}), we have
\[ \frac{2u_1^2}{1+\sqrt 2} = x_1-x_2 = \frac{2u_2^2}{-1+\sqrt 2}, \ \ \frac{2u_3^2}{1+\sqrt 2} = x_1+x_2 = \frac{2u_4^2}{-1+\sqrt 2}.\]
Since $2u_1u_2 = g_2 = x_2-x_1\leq 0$ and $2u_3u_4 = g_4 = x_1+x_2\geq 0$, it follows that
\begin{equation}\label{uvariables} u_2 = -(\sqrt 2-1)u_1,  \ \ u_4  =  (\sqrt 2-1)u_3.\end{equation}
From the second canonical transformation (\ref{secondcantrans}), we have
\begin{align*}
v_1 & = \sqrt 2(\omega_1-\omega_2) u_1, & v_2 & = -(2-\sqrt 2)(\omega_1-\omega_2)u_1, \\
v_3 & = \sqrt 2(\omega_1+\omega_2) u_3, & v_4 & = (2-\sqrt 2)(\omega_1+\omega_2)u_3.
\end{align*}
These imply that
\begin{equation}\label{vvariables} v_2 = -(\sqrt 2 -1)v_1,\ \ v_4  = (\sqrt 2-1)v_3.\end{equation}
Substitution into the Hamiltonian system of equations with Hamiltonian $\hat \Gamma$ (and with $m=1$) gives
\begin{align*} u_1^\prime = & \frac{4-2\sqrt 2}{4} v_1u_3^2,
& u_2^\prime = & -\frac{(\sqrt 2 - 1)(4-2\sqrt 2)}{4} v_1u_3^2, \\
u_3^\prime = & \frac{4-2\sqrt 2}{4}v_3u_1^2, 
& u_4^\prime = & \frac{(\sqrt 2-1)(4-2\sqrt 2)}{4} v_3u_1^2,
\end{align*}
$\hat E^\prime = 0$, and
\begin{align*}
v_1^\prime = & -\frac{(4-2\sqrt 2)u_1v_3^2}{4}  + 4u_1+ \frac{4u_1u_3^2}{\sqrt{u_1^4+u_3^4}} - \frac{4u_1^5u_3^2}{\big(u_1^4+u_3^4\big)^{3/2}} + 2(4-2\sqrt 2)\hat E u_1u_3^2, \\
v_2^\prime = & -(\sqrt 2-1)\Bigg[ -\frac{(4-2\sqrt 2)u_1v_3^2}{4}  + 4u_1+ \frac{4u_1u_3^2}{\sqrt{u_1^4+u_3^4}} - \frac{4u_1^5u_3^2}{\big(u_1^4+u_3^4\big)^{3/2}}  \\
& \hspace{0.8in} + 2(4-2\sqrt 2)\hat E u_1u_3^2\Bigg], \\
v_3^\prime = &  -\frac{(4-2\sqrt 2)u_3v_1^2}{4}  + 4u_3+ \frac{4u_1^2u_3}{\sqrt{u_1^4+u_3^4}} - \frac{4u_1^2u_3^5}{\big(u_1^4+u_3^4\big)^{3/2}} + 2(4-2\sqrt 2)\hat E u_1^2u_3, \\
v_4^\prime  = &\ (\sqrt 2-1)\Bigg[-\frac{(4-2\sqrt 2)u_3v_1^2}{4}  + 4u_3+ \frac{4u_1^2u_3}{\sqrt{u_1^4+u_3^4}} - \frac{4u_1^2u_3^5}{\big(u_1^4+u_3^4\big)^{3/2}} \\
& \hspace{0.7in} + 2(4-2\sqrt 2)\hat E u_1^2u_3\Bigg], \\
t^\prime  = &\ (4-2\sqrt 2)^2u_1^2u_3^2.
\end{align*}
Because of Equations (\ref{uvariables}) and (\ref{vvariables}), the equations in $u_2^\prime$, $u_4^\prime$, $v_2^\prime$, and $v_4^\prime$ duplicate those in $u_1^\prime$, $u_3^\prime$, $v_1^\prime$, and $v_3^\prime$.
The Ans\"atz $x_4=x_1$, $x_3=x_2$ with $\vert x_2\vert\leq x_1$, therefore leads to the reduced system of equations,
\begin{align*}
u_1^\prime = & \frac{4-2\sqrt 2}{4} v_1u_3^2, \\
u_3^\prime = & \frac{4-2\sqrt 2}{4}v_3u_1^2, \\
\hat E^\prime = & 0, \\
v_1^\prime = & -\frac{(4-2\sqrt 2)u_1v_3^2}{4}  + 4u_1+ \frac{4u_1u_3^2}{\sqrt{u_1^4+u_3^4}} - \frac{4u_1^5u_3^2}{\big(u_1^4+u_3^4\big)^{3/2}} + 2(4-2\sqrt 2)\hat E u_1u_3^2, \\
v_3^\prime = &  -\frac{(4-2\sqrt 2)u_3v_1^2}{4}  + 4u_3+ \frac{4u_1^2u_3}{\sqrt{u_1^4+u_3^4}} - \frac{4u_1^2u_3^5}{\big(u_1^4+u_3^4\big)^{3/2}} + 2(4-2\sqrt 2)\hat E u_1^2u_3,\\
t^\prime = & (4-2\sqrt 2)^2u_1^2u_3^2.
\end{align*}
Scale the value of $\hat E$ by
\begin{equation}\label{energyEtilde} \tilde E = \frac{\hat E}{4-2\sqrt 2},\end{equation}
and define
\[ \tilde \Gamma = \frac{4-2\sqrt 2}{8}\big( v_1^2u_3^2+ v_3^2u_1^2\big) - 2(u_1^2+u_3^2) - \frac{2u_1^2u_3^2}{\sqrt{u_1^4+u_3^4}} - (4-2\sqrt 2)^2\tilde E u_1^2u_3^2.\]
It is straight-forward to check that the reduced system of equations satisfies
\[ u_i^\prime = \frac{\partial\tilde\Gamma}{\partial v_i}, \ \ v_i^\prime = -\frac{\partial\tilde\Gamma}{\partial u_i},\ \ i=1,2,\ \  \tilde E^\prime = \frac{\partial\tilde\Gamma}{\partial t}, \ \ t^\prime = - \frac{\partial\tilde\Gamma}{\partial \tilde E}.\]
Thus the system of reduced equations is Hamiltonian.

We will simplify the Hamiltonian $\tilde \Gamma$ by a linear symplectic transformation with a multiplier $\mu \ne 1$. Define new coordinates $(Q_1,Q_2,E,P_1,P_2,\tau)$ by
\begin{align}\label{lstmu1} u_1 & = \frac{Q_1}{2^{1/4}}, & v_1 & = \frac{P_1}{2\sqrt{\sqrt 2 - 1}}, \\
\label{lstmu2} u_3 & = \frac{Q_2}{2^{1/4}}, & v_3 & = \frac{P_2}{2\sqrt{\sqrt 2 - 1}}, \\
\label{lstmu3} \tilde E & = \frac{2E}{(4-2\sqrt 2)^2}, & t & = 2^{3/4}(\sqrt 2 - 1)^{3/2}\ \tau.
\end{align}
This is a linear symplectic change of coordinates with multiplier
\begin{equation}\label{mu} \mu=\frac{1}{2^{5/4}\sqrt{\sqrt 2-1}}.\end{equation}
Under this linear symplectic transformation and the accompanying scaling $\sigma= s/\mu$ of the independent variable $s$, the Hamiltonian $\tilde\Gamma$ becomes
\[ \Gamma = \frac{1}{16}\big(P_1^2Q_2^2 + P_2^2Q_1^2\big) -\sqrt2(Q_1^2+Q_2^2) - \frac{\sqrt 2 Q_1^2Q_2^2}{\sqrt{Q_1^4+Q_2^4}} - EQ_1^2Q_2^2.\]
The reduced system of equations is the Hamiltonian system of equations with Hamiltonian $\Gamma$,
\begin{equation}\label{gamma} \frac{dQ_i}{d\sigma} = \frac{\partial \Gamma}{\partial P_i}, \ \ \frac{dP_i}{d\sigma}=-\frac{\partial \Gamma}{\partial Q_i}, \ \ i=1,2, \ \ \frac{dE}{d\sigma}  = 0, \ \ \frac{d\tau}{d\sigma}  = Q_1^2Q_2^2.\end{equation}

The function $\Gamma$ is a regularized Hamiltonian for the fully symmetric planar four-body equal mass problem with the bodies located at $(x_1,x_2)$, $(x_2,x_1)$, $(-x_1,-x_2)$, and $(-x_2,-x_1)$ (see \cite{OY3}). On the level set $\Gamma=0$, the solutions have energy $E$. One regularized simultaneous binary collision occurs when $Q_1=0$ and $Q_2\ne 0$, for which $\Gamma=0$ implies $P_1^2 = 16\sqrt 2$, and for which the transformation between $Q_1,Q_2$ and $x_1,x_2$ implies $x_1-x_2=0$ and $x_1+x_2\ne 0$. The other regularized simultaneous binary collision occurs when $Q_1\ne 0$ and $Q_2=0$, for which $\Gamma=0$ implies $P_2^2 = 16\sqrt 2$, and for which the transformation between $Q_1,Q_2$ and $x_1,x_2$ implies $x_1-x_2\ne 0$ and $x_1+x_2=0$. Total collapse occurs when $Q_1=0$ and $Q_2=0$, and is the only singularity in $\Gamma$ that is not regularized. A solution $Q_1(\sigma)$, $Q_2(\sigma)$, $P_1(\sigma)$, $P_2(\sigma)$, $\sigma\in[0,\sigma_0]$, for $\sigma_0>0$, of the Hamiltonian system equations with Hamiltonian $\Gamma$ is {\it nonsingular}\, if it avoids total collapse, i.e., $Q_1^4+Q_2^4\ne 0$ for all $\sigma\in[0,\sigma_0]$.

The following result is from \cite{OY3}. The proof of it is a consequence of four equal mass bodies starting at $(x_1,x_2)$, $(x_2,x_1)$, $(-x_1,-x_2)$, and $(-x_2,-x_1)$ with $x_1=1$ and $x_2=0$, and with the momenta $(0,\vartheta)$, $(\vartheta,0)$, $(0,-\vartheta)$, and $(-\vartheta,0)$ for any $\vartheta>0$, always having a simultaneous binary collision on the line $x_2=x_1$ at a time $t_0>0$ continuously depending on $\vartheta$, such that the cluster velocity $\dot x_1(t_0)+\dot x_2(t_0)$ is a continuous function of $\vartheta$.

\begin{lemma}\label{reducedexistence} There exists $\vartheta>0$, $\sigma_0>0$, and a nonsingular solution $Q_1(\sigma)$, $Q_2(\sigma)$, $P_1(\sigma)$, $P_2(\sigma)$, $\sigma\in[0,\sigma_0]$, of the Hamiltonian system of equations with Hamiltonian $\Gamma$ on the level set $\Gamma=0$ such that
\begin{equation}\label{initial0} Q_1(0) = 1, \ \ Q_2(0) = 1, \ \ P_1(0) = -\vartheta,\ \ P_2(0) = \vartheta,\end{equation}
\begin{equation}\label{energyE} E = \frac{\vartheta^2-16\sqrt 2 - 8}{8}<0,\ \ \tau(\sigma) = \int_0^\sigma Q_1^2(y)Q_2^2(y)\ dy,\end{equation}
and 
\begin{equation}\label{terminal} Q_1(\sigma_0) = 0, \ \ Q_2(\sigma_0) >0, \ \ P_1(\sigma_0) =-4 (2^{1/4}), \ \ P_2(\sigma_0)=0.\end{equation}
\end{lemma}

This Lemma gives the existence of a solution of a boundary value problem for the Hamiltonian system of equations with Hamiltonian $\Gamma$. It is this solution whose symmetric extension gives a symmetric periodic SBC orbit in the PPS4BP with equal masses.

\begin{theorem}\label{existence} Fix $m=1$. For each $\hat E<0$, there exists a time-reversible periodic regularized SBC orbit
\[ \gamma(s) =(u_1(s),u_2(s),u_3(s),u_4(s),v_1(s),v_2(s),v_3(s),v_4(s))\]
with period $T>0$, angular momentum $A=0$, and a symmetry group isomorphic to $D_4$, for the Hamiltonian system of equations with Hamiltonian $\hat\Gamma$ on the level set $\hat\Gamma=0$ such that distinct regularized SBCs occur at $s=T/8,3T/8,5T/8,7T/8$. This periodic orbit corresponds to a symmetric periodic singular orbit
\[ (x_1(t),\ x_2(t),\ x_3(t),\ x_4(t),\ \omega_1(t),\ \omega_2(t),\ \omega_3(t),\ \omega_4(t))\]
with energy $\hat E$ for the PPS4BP with equal masses $m=1$ where for all $t$,
\[ x_4(t) = x_1(t),\ x_3(t) = x_2(t),\ \vert x_2(t)\vert\leq x_1(t),\ \omega_4(t)=\omega_1(t),\ \omega_3(t)=\omega_2(t),\]
with initial conditions
\[x_1(0)>0,\ x_2(0)=0,\ \omega_1(0)=0,\ \omega_2(0)>0,\]
and period 
\[ R=\int_0^{T/2} \big(u_1^2(s)+u_2^2(s)\big)\big(u_3^2(s)+u_4^2(s)\big)\ ds,\]
where for $t\in[0,R]$, the only singularities are two distinct SBC's occurring at $t=R/4$ and $t=3R/4$.
\end{theorem}

\begin{proof} By Lemma \ref{reducedexistence}, let $Q_1(\sigma)$, $Q_2(\sigma)$, $P_1(\sigma)$, $P_2(\sigma)$, $\sigma\in[0,\sigma_0]$, be the nonsingular solution of the Hamiltonian system of equations (\ref{gamma}) with Hamiltonian $\Gamma$ on the level set $\Gamma=0$, whose properties are given in (\ref{initial0}), (\ref{energyE}), and (\ref{terminal}). Using the scaling $\sigma =s/\mu$, set $s_0 = \mu \sigma_0$. By the linear symplectic transformation given in (\ref{lstmu1}), (\ref{lstmu2}), and (\ref{lstmu3}) with multiplier $\mu$ (as given in (\ref{mu})), we have
\[ u_1(s) = \frac{Q_1(s/\mu)}{2^{1/4}}, \ \ u_3(s) = \frac{Q_2(s/\mu)}{2^{1/4}}, \ \ v_1(s) = \frac{P_1(s/\mu)}{2\sqrt{\sqrt 2 -1}}, \ \ v_3(s) = \frac{P_2(s/\mu)}{2\sqrt{\sqrt 2 -1}}.\]
By Equations (\ref{uvariables}) and (\ref{vvariables}), we have
\[ u_2(s) = \frac{-(\sqrt 2-1)Q_1(s/\mu)}{2^{1/4}}, \ \ u_4(s) = \frac{(\sqrt 2-1)Q_2(s/\mu)}{2^{1/4}}, \]
\[ v_2(s) = \frac{-\sqrt{\sqrt 2-1} P_1(s/\mu)}{2}, \ \ v_4(s) = \frac{\sqrt{\sqrt 2 -1} P_2(s/\mu)}{2}.\]
Set $\gamma(s) = (u_1(s),u_2(s),u_3(s),u_4(s),v_1(s),v_2(s),v_3(s),v_4(s))$, $s\in[0,s_0]$. From the Equations (\ref{energyEtilde}), (\ref{lstmu3}), (\ref{energyE}) in $E$, $\hat E$, $\tilde E$, and $\vartheta$, we obtain
\[ \hat E = \frac{(2+\sqrt 2)\vartheta^2}{16} - 3 - \frac{5\sqrt 2}{2}<0.\]
With this value of $\hat E$, it follows that the value of $\hat\Gamma$ at $\gamma(0)$ is $0$. Set
\[ \zeta_1 = \frac{1}{2^{1/4}}>0, \ \ \zeta_2 = \frac{-(\sqrt 2-1)}{2^{1/4}}<0,\]
\[ \rho_1 = \frac{-\vartheta}{2\sqrt{\sqrt 2 - 1}}<0, \ \ \rho_2 = \frac{\vartheta\sqrt{\sqrt 2-1} }{2}>0.\]
With $Q_1(\sigma)$, $Q_2(\sigma)$, $P_1(\sigma)$, $P_2(\sigma)$, $\sigma\in[0,\sigma_0]$, being a nonsingular solution of the Hamiltonian system of equations with Hamiltonian $\Gamma$, we have $Q_1^4(\sigma)+Q_2^4(\sigma)\ne 0$ for all $\sigma\in[0,\sigma_0]$. From this it follows for all $s\in[0,s_0]$ that
\[ M_5 = u_1^2(s) - u_2^2(s) + u_3^2(s) - u_4^2(s) = (2-\sqrt 2)\big[Q_1^2(s/\mu) + Q_2^2(s/\mu)\big] \ne 0,\]
and
\[ M_8 = 2u_1(s)u_2(s) - 2u_3(s)u_4(s) = - (2-\sqrt 2)\big[ Q_1^2(s/\mu) + Q_2^2(s/\mu)\big] \ne 0.\]
These imply that $M_5^2+M_6^2\ne 0$ and $M_7^2+M_8^2\ne 0$ for all $s\in[0,s_0]$. Thus the function $\gamma(s)$ is a nonsingular solution of the Hamiltonian system of equations (\ref{regularized}) with Hamiltonian $\hat\Gamma$ on the level set $\hat\Gamma=0$ that satisfies
\[ u_1(0) = \zeta_1,\ \ u_2(0) = \zeta_2 \ \ u_3(0) = \zeta_1, \ \ u_4(0) = -\zeta_2,\]
\[ v_1(0) = \rho_1, \ \ v_2(0) = \rho_2, \ \ v_3(0) = -\rho_1, \ \ v_4(0) = \rho_2,\]
\[ u_1(s_0) = 0, \ \ u_2(s_0) = 0, \ \ u_3(s_0) > 0, \ \ u_4(s_0) > 0,\]
\[ v_1(s_0) <0, \ \ v_2(s_0) >0, \ \ v_3(s_0) = 0, \ \ v_4(s_0) = 0. \vspace{0.05in}\]
By Lemma \ref{symmetry}, the solution $\gamma(s)$ extends to a $T=8s_0$ periodic solution, call it $\gamma(s)$, with a $D_4$ symmetry group generated by the symmetries $S_F$ and $S_G$, and four distinct regularized simultaneous binary collisions at $s=s_0,3s_0,5s_0,7s_0$, for which
\begin{align*} & u_1^2(s_0)+u_2^2(s_0)=0,\ u_3^2(s_0)+u_4^2(s_0)\ne 0, \\
& u_1^2(3s_0)+u_2^2(3s_0)\ne 0,\ u_3^2(3s_0)+u_4^2(3s_0)= 0, \\
& u_1^2(5s_0)+u_2^2(5s_0)=0,\ u_3^2(5s_0)+u_4^2(5s_0)\ne 0, \\
& u_1^2(7s_0)+u_2^2(7s_0)\ne 0,\ u_3^2(7s_0)+u_4^2(7s_0)=0.
\end{align*}
Since $Q_1^4(\sigma)+Q_2^4(\sigma)\ne 0$ for all $\sigma\in[0,\sigma_0]$, it follows that $(u_1^2(s)+u_2^2(s))(u_3^2(s)+u_4^2(s))\ne 0$ for $s\in[0,T]$ except $s=s_0,3s_0,5s_0,7s_0$. The regularizing change of time (\ref{timechange}),
\[ \frac{dt}{ds} = (u_1^2(s)+u_2^2(s))(u_3^2(s)+u_4^2(s)),\]
defines $t$ as an invertible differentiable function of $s$, i.e., $t=\theta(s)$ with $\theta(0)=0$ and $\theta^\prime(s)=0$ when $s=(2k+1)s_0$ for $k\in{\mathbb Z}$. The symmetry $S_F$ satisfies $S_F\gamma(s) = \gamma(s+2s_0)$ and $-\gamma(s)=S_F^2\gamma(s) = \gamma(s+4s_0)$. The symmetry $S_G$ satisfies $S_G\gamma(s) = \gamma(2s_0-s)$, and so $\gamma(s)$ has a time-reversing symmetry. The angular momentum (\ref{ram}) for $\gamma(s)$ at $s=0$ is
\[ A =\frac{1}{2}\big[ -v_1u_2 + v_2u_1 - v_3u_4 + v_4u_3\big] = \rho_2\zeta_1-\rho_1\zeta_2 = 0.\]

The extended solution $\gamma(s)$ gives a singular symmetric solution
\[ z(t)=(x_1(t),x_2(t),x_3(t),x_4(t),\omega_1(t),\omega_2(t),\omega_3(t),\omega_4(t)),\] of the PPS4BP with $m=1$. Under the Ans\"atz, the components of $z(t)$ satisfies $x_4(t) =x_1(t)$, $x_3(t)=x_2(t)$, $\vert x_2(t)\vert\leq x_1(t)$, $\omega_4(t)=\omega_1(t)$, $\omega_3(t) = \omega_2(t)$, where
\begin{align*} x_1(t) & = \frac{u_1^2(s) - u_2^2(s) + u_3^2(s) - u_4^2(s)}{2} \\
x_2(t)  & = u_1(s)u_2(s) + u_3(s)u_4(s), \\
\omega_1(t) & = \frac{v_1(s)u_1(s) - v_2(s)u_2(s)}{2(u_1^2(s)+u_2^2(s))} + \frac{v_3(s)u_3(s) - v_4(s)u_4(s)}{2(u_3^2(s)+u_4^2(s))},\\
\omega_2(t) & = \frac{ v_1(s)u_2(s) + v_2(s)u_1(s)}{2(u_1^2(s)+u_2^2(s))} + \frac{v_3(s)u_4(s)+ v_4(s)u_3(s)}{2(u_3^2(s)+u_4^2(s))},
\end{align*}
for $s=\theta^{-1}(t)$. The components of the extended solution $\gamma(s)$ satisfy $u_3(0)= u_1(0)$, $u_4(0)=-u_2(0)$, $u_1(0)u_2(0)<0$, $\vert u_2(0)\vert = (\sqrt 2 -1)\vert u_1(0)\vert$, $v_3(0) = -v_1(0)$, $v_4(0) = v_2(0)$, 
\[ v_1(0)u_2(0) + v_2(0)u_1(0) = \rho_1\zeta_2+\rho_2\zeta_1 = \frac{ \sqrt{\sqrt 2 -1}\ \vartheta}{2^{1/4}}>0,\]
and
\[ v_2(0)u_2(0) - v_1(0)u_1(0) = \rho_2\zeta_2 - \rho_1\zeta_1 =  \frac{ \sqrt{\sqrt 2 -1}\ \vartheta}{2^{1/4}}.\]
From Lemma \ref{initial}, it follows that $x_1(0)>0$, $x_2(0)=0$, $\omega_1(0)=0$, and $\omega_2(0)>0$. Set $R=\theta(T/2)$. Since $\gamma(4s_0) = -\gamma(0)$, it follows that $x_1(R) = x_1(0)$, $x_2(R)=x_2(0)$, $\omega_1(R)=\omega_1(0)$, and $\omega_2(R)=\omega_2(0)$. Thus the singular symmetric solution $z(t)$ has period $R$. By the construction of the extension of $\gamma(s)$ given in Lemma \ref{symmetry}, there holds
\[ \int_{ks_0}^{(k+1)s_0} (u_1^2(s)+u_2^2(s))(u_3^2(s)+u_4^2(s))\ ds = \int_0^{s_0} (u_1^2(s)+u_2^2(s))(u_3^2(s)+u_4^2(s))\ ds\]
for all $k=1,\dots,7$. This implies that $R/4 = \theta((k+1)s_0) -\theta(ks_0)$ for all $k=0,1,\dots,7$. The first regularized SBC for $\gamma(s)$ occurs at $s=s_0$, and this corresponds to $t=\theta(s_0)=R/4$. The next regularized SBC for $\gamma(s)$ occurs at $s=3s_0$, and this corresponds to
\[ t = \theta(3s_0) = (\theta(3s_0)-\theta(2s_0)) + (\theta(2s_0) -\theta(s_0)) + \theta(s_0) = \frac{R}{4}+\frac{R}{4} + \frac{R}{4} = \frac{3R}{4}.\]
Similarly, the regularized SBCs for $\gamma(s)$ occurring at $s=5s_0,7s_0$ correspond to $t=5R/4, 7R/4$. Hence, for $t\in[0,R]$, the periodic solution $z(t)$ has SBCs as its only singularities, and these occur at $t=R/4,3R/4$.

For a fixed but arbitrary $\epsilon>0$, the value of $\epsilon^{-2}\hat E$ is a fixed but arbitrary negative real number. By Lemma \ref{energy}, the scaled extended solution
\[ \gamma_\epsilon(s) = (\epsilon u_1(\epsilon s),\epsilon u_2(\epsilon s),\epsilon u_3(\epsilon s),\epsilon u_4(\epsilon s),v_1(\epsilon s),v_2(\epsilon s),v_3(\epsilon s),v_4(\epsilon s))  \]
is a periodic solution for the Hamiltonian system of equations (\ref{regularized}) with Hamiltonian $\hat\Gamma$ on the level set $\hat\Gamma=0$, having period $\epsilon^{-1}T$ and energy $\epsilon^{-2}\hat E<0$. By an argument similar to above, $\gamma_\epsilon(s)$ satisfies the required conditions.
\end{proof}

\section{Numerical Estimates in the Equal Mass Case}
In the equal mass case, there is by Theorem \ref{existence} a time-reversible periodic orbit $\gamma(s)$ for the Hamiltonian system of equations with Hamiltonian $\hat\Gamma$ on the level set $\hat\Gamma=0$ with period $T$. The components $u_1(s)$, $u_2(s)$, $u_3(s)$, $u_4(s)$, $v_1(s)$, $v_2(s)$, $v_3(s)$, $v_4(s)$ of $\gamma(s)$ satisfy Equations (\ref{uvariables}) and (\ref{vvariables}). The $D_4$ symmetry group of $\gamma(s)$ is generated by $S_F\gamma(s) = \gamma(s+T/4)$ and $S_G\gamma(s) = \gamma(T/4 -s)$.  Under the linear symplectic transformation (\ref{lstmu1}), (\ref{lstmu2}), and (\ref{lstmu3}) with multiplier $\mu$ (as given by (\ref{mu})), this gives a periodic orbit $Q_1(\sigma)$, $Q_2(\sigma)$, $P_1(\sigma)$, $P_2(\sigma)$ of the Hamiltonian system of equations (\ref{gamma}) with Hamiltonian $\Gamma$ on the level set $\Gamma=0$, which by Lemma \ref{reducedexistence} satisfies $Q_1(0)=1$, $Q_2(0)=1$, $P_1(0) = -\vartheta$, $P_2(0)=\vartheta$, $Q_1(\sigma_0)=0$, $Q_2(\sigma_0)>0$, $P_1(\sigma_0)=-4(2^{1/4})$, and $P_2(\sigma_0)=0$ for some $\sigma_0>0$ and  $\vartheta>0$. In \cite{BOYS}, we numerically estimated
\[ \sigma_0 = 1.62047369909693,\ \ \vartheta = 2.57486992651942.\]
The period of this periodic orbit is $8\sigma_0\approx 12.96378959$ and its energy is $E\approx -2.999682732$. From the linear symplectic transformation (\ref{lstmu1}), (\ref{lstmu2}), and (\ref{lstmu3}), with multiplier $\mu$ and from Equations (\ref{uvariables}) and (\ref{vvariables}), we have for the values of components of $\gamma(0)$ the exact
\[ u_1(0) = u_3(0)= 2^{-1/4} ,\ \ -u_2(0) = u_4(0) = (\sqrt 2 -1)2^{-1/4},\]
and the estimates
\begin{align*} v_1(0) & = -v_3(0) = -\frac{\vartheta}{2\sqrt{\sqrt 2-1}}\approx-2.000382939,\\
v_2(0) & = v_4(0) = \frac{\vartheta\sqrt{\sqrt 2-1}}{2} \approx 0.8285857433.
\end{align*}
Since $\sigma=s/\mu$, the period of $\gamma(s)$ is $T\approx 8.469003682$. From Equation (\ref{energyEtilde}), the value of the energy for $\gamma(s)$ is $\hat E\approx -5.120778733$.

\begin{figure}\scalebox{0.38}{\includegraphics[trim=20mm 0mm 0mm 0mm,clip]{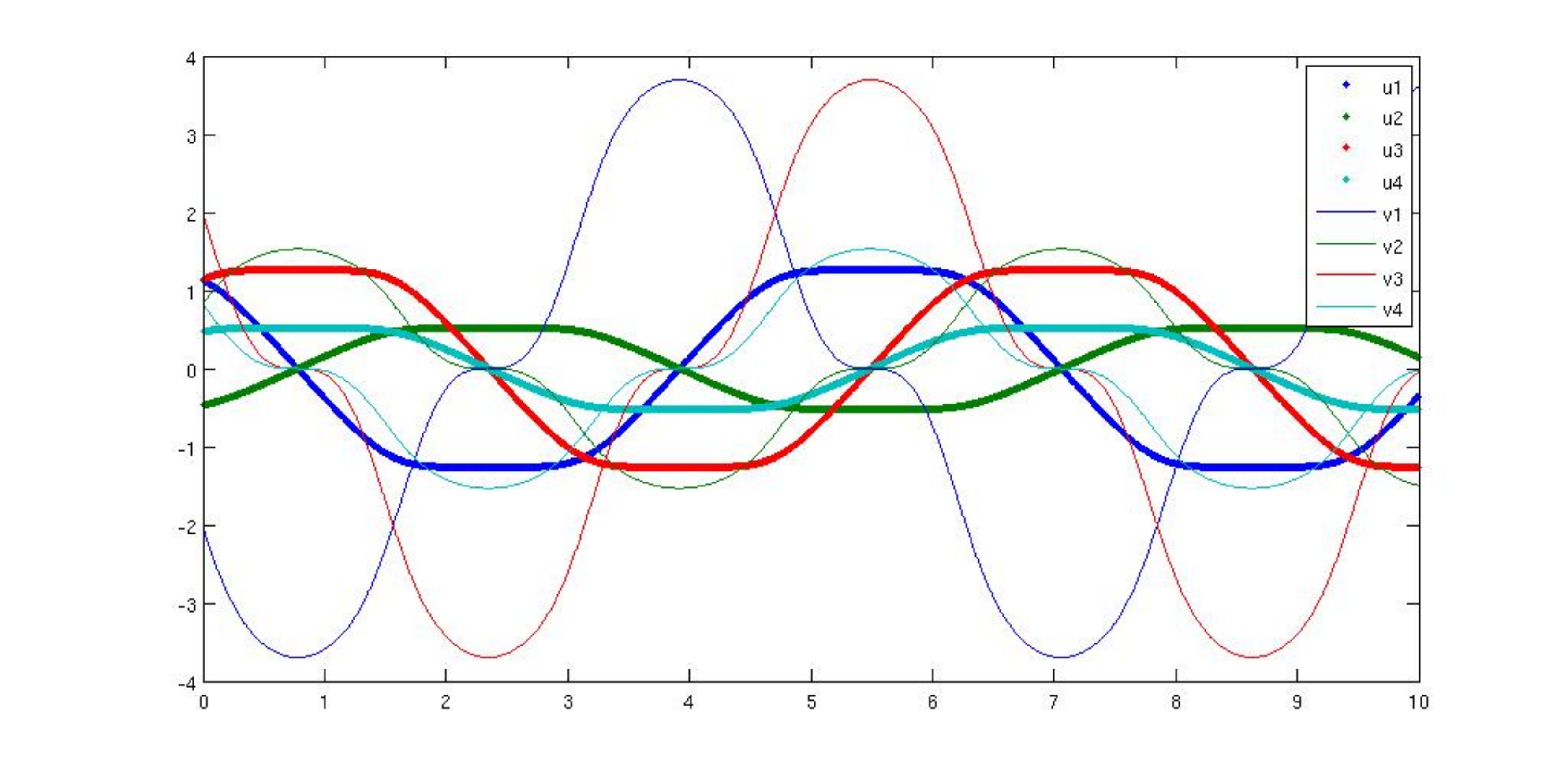}}
\caption{The $u_1$, $u_2$, $u_3$, $u_4$, $v_1$, $v_2$, $v_3$, $v_4$ coordinates of the symmetric periodic SBC orbit $\gamma_\epsilon(s)$ for $m=1$ and an $\epsilon>1$.}\label{figure2}
\end{figure}

The components of the scaled periodic orbit $\gamma_\epsilon(s)$ for an $\epsilon>1$, shown in Figure \ref{figure2}, satisfy the symmetries $S_F\gamma_\epsilon(s) = \gamma_\epsilon(s+T_\epsilon/4)$ and $S_G\gamma_\epsilon(s) = \gamma_\epsilon(T_\epsilon/4 -s)$. We choose $\epsilon$ so that $T_\epsilon=2\pi$, and check the linear stability of $\gamma_\epsilon(s)$ (see Theorem \ref{stability}). Using a Runge-Kutta order $4$ algorithm with a fixed time step of $2\pi/50000$, we computed $X_\epsilon(2\pi)$. Two of the eigenvalues of $X_\epsilon(2\pi)$ are $1$ by Theorem \ref{eigenvalueone}. Numerical estimates of the remaining eigenvalues of $X_\epsilon(2\pi)$ are
\begin{align*}  -0.9888731375 & \pm 0.1487612779i, \\  -0.9973584383 & \pm 0.07263708002i, \\  0.9999060579 & \pm 0.01370676220i,
\end{align*}
that have modulus one. Thus numerically, the periodic orbits $\gamma_\epsilon(s)$ are linearly stable for all $\epsilon>0$. The first complex conjugate pair of eigenvalues for $X_\epsilon(2\pi)$ matches the complex conjugate pair of characteristic multipliers for the periodic orbit $Q_1(\sigma)$, $Q_2(\sigma)$, $P_1(\sigma)$, $P_2(\sigma)$ of the Hamiltonian system of equations (\ref{gamma}) with Hamiltonian $\Gamma$, corresponding to $\gamma_1(s)=\gamma(s)$, where we computed  \cite{BOYS} the real part of the complex conjugate pair of modulus one to be $-0.9888840619$.  Because, by a lengthly computation, $J\nabla^2\hat\Gamma(\gamma_\epsilon(s))$ is not block diagonal, the last two complex conjugate pairs of eigenvalues of $X_\epsilon(2\pi)$ are not repeats of the characteristic multipliers of the periodic orbit $Q_1(\sigma)$, $Q_2(\sigma)$, $P_1(\sigma)$, $P_2(\sigma)$.

Figure \ref{figure1} illustrates the curves in the physical plane that the four equal masses follow in the linearly stable SBC orbit $z(t)=(x_1(t)$, $x_2(t)$, $x_3(t)$, $x_4(t)$, $\omega(t)$, $\omega_2(t)$, $\omega_3(t)$, $\omega_4(t))$ of the PPS4BP, corresponding to $\gamma_\epsilon(s)$ with $\epsilon=1/(2-\sqrt 2)$. The initial conditions for $z(t)$ are
\begin{align*} & x_1(0) = x_4(0) = 1, \ \ x_2(0) = x_3(0) = 0, \\
& \omega_1(0) = \omega_4(0) = 0, \ \ \omega_2(0) = \omega_3(0) \approx 1.287434964.
\end{align*}
The value of $\epsilon=1/\sqrt{2-\sqrt 2}$ here for the scaling is determined by the equation $x_1(0)= \epsilon^2(u_1^2(0)-u_2^2(0))$ coming from  Lemma \ref{initial} and the canonical transformations (\ref{firstcantrans}) and (\ref{secondcantrans}) applied to the scaled periodic solution $\gamma_\epsilon(s)$, together with the initial condition $x_1(0)=1$, where $t=0$ corresponds to $s=0$. The value of the Hamiltonian $H$ along $z(t)$ is $\epsilon^{-2}\hat E\approx -2.999682732$.

\section{Numerical Estimates in the Unequal Mass Case}

We numerically continue to $0<m<1$ the time-reversible periodic regularized simultaneous binary collision orbit for the Hamiltonian system of equations with Hamiltonian $\hat\Gamma$ on the level set $\hat\Gamma=0$. The analytic existence of this orbit is given by Theorem \ref{existence} for $m=1$. We also numerically estimate the monodromy matrices and their eigenvalues for these continued time-reversible periodic regularized simultaneous binary collision orbits for $0< m<1$. For $m=1$, we assume by Lemma \ref{energy} that the time-reversible periodic regularized simultaneous binary collision orbit $\gamma(s;1)$ with the $D_4$ symmetry has period $T=2\pi$ and energy $\hat E\approx-2.818584789$. Recall that the $D_4$ symmetry group of $\gamma(s;1)$ is generated by $S_F\gamma(s;1) = \gamma(s+\pi/2;1)$ and $S_G\gamma(s;1) = \gamma(\pi/2-s;1)$. We assume without loss of generality by Lemma \ref{energy} and Lemma \ref{symmetry}, that the continued time-reversible periodic regularized simultaneous binary collision orbit $\gamma(s;m)$ for $0<m<1$, has period $T=2\pi$, energy $\hat E(m)$, and a $D_4$ symmetry group generated by $S_F\gamma(s;m)=\gamma(s+\pi/2;m)$ and $S_G\gamma(s;m) = \gamma(\pi/2-s;m)$. We assume $\hat E(m)$ depends continuously on $0<m\leq 1$, with $\hat E(1)=\hat E$. We shift the regularized time variable $s$ to $s+\pi/4$, so that for
\begin{align*} \gamma(s;m) = & (u_1(s;m),u_2(s;m),u_3(s;m),u_4(s;m), \\ & v_1(s;m),v_2(s;m),v_3(s;m),v_4(s;m)),\end{align*}
the value $s=0$ now corresponds to the first simultaneous binary collision, i.e., $u_1(0;m)=0$, $u_2(0;m)=0$, $u_3(0;m)\ne0$, $u_4(0;m)\ne 0$, $v_1(0;m) \ne 0$, $v_2(0;m) \ne 0$, $v_3(0;m) =0$, and $v_4(0,m) = 0$. We approximate the continued time-reversible periodic orbits $\gamma(s;m)$, $0<m\leq 1$, by the trigonometric polynomials,
\begin{align*} & u_1(s;m) = \sum_{i=1}^n a_i \sin((2i-1)s), \ \ u_2(s;m) =  -\sum_{i=1}^n b_i \sin((2i-1)s),\\ & u_3(s;m) = u_1(s-\pi/2;m), \ \ u_4(s;m) = u_2(s+\pi/2;m), \\
& v_1(s;m) =  \sum_{i=1}^n c_i \cos((2i-1)s), \ \ v_2(s;m) = -\sum_{i=1}^n d_i \cos((2i-1)s),\\
& v_3(s;m) = v_1(s-\pi/2;m), \ \ v_4(s;m) = v_2(s+\pi/2;m),
\end{align*}
for a positive integer $n$ and constants $a_i$, $b_i$, $c_i$, $d_i$, $i=1,\dots,n$, that are assumed to be continuous functions of $m$. The presence of odd positive integer frequencies $2i-1$ in these trigonometric polynomials is to ensure that the period functions defined by them have the $D_4$ symmetry group generated by $S_F$ and $S_G$. So in particular, the periodic orbits in the continuation are time-reversible for all $0<m<1$. An numerical estimate of the periodic solution $\gamma(s;m)$ is found through the variational approach of minimizing the functional
\[ L=\int_0^{2\pi} \Vert \gamma^{\, \prime}(s;m) - J\nabla\hat\Gamma(\gamma(s;m))\Vert\ ds\]
over the space $\mathbb{R}^{4n}$ of coefficients $a_i$, $b_i$, $c_i$, $d_i$, $i=1,\dots,n$, for an appropriate choice of $n$. The use of trigonometric polynomials for numerically approximating periodic solutions is a classic approach (see, for example, Sim\'o \cite{Si}).

The numerical algorithm for finding a trigonometric polynomial that approximates the periodic solution $\gamma(s;m)$, $0<m<1$, proceeds in two steps.  First, we consider a guess for the set of values for $a_i$, $b_i$, $c_i$, $d_i$, $i=1,\dots,n$, as well as a guess for $\hat{E}(m)$.  Starting with a reasonably low number of terms, $n = 5$, we let a numerical minimization algorithm (in this case, MATLAB's \verb|fminunc|) find the minimum of $L$ near the starting guess.  Then we add an additional non-zero term to each of the trigonometric polynomials, and the minimizing solution from the previous iteration is used as a starting guess for the next iteration.  This process continues until we reach $n = 10$. Second, since the Hamiltonian equations with Hamiltonian $\hat\Gamma$ requires a specified value of $\hat E(m)$ to compute $L$, we need to make certain that we are getting a good estimate of $\hat E(m)$. We evaluate $\hat{\Gamma}$ at the point $\gamma(\pi/4;m)$, away from simultaneous binary collisions. If this value is not sufficiently close to $0$ (within about $5\times 10^{-10}$ of $0$), we adjust $\hat E(m)$ using the bisection method in a small interval about the initial guess of $\hat E(m)$ until $\hat{\Gamma}(\gamma(\pi/4;m))$ is sufficiently close to $0$, re-minimizing the trigonometric polynomial approximation of $\gamma(s;m)$ for each new choice of $\hat E(m)$.

\begin{figure}\scalebox{0.4}{\includegraphics[trim=15mm 0mm 0mm 0mm]{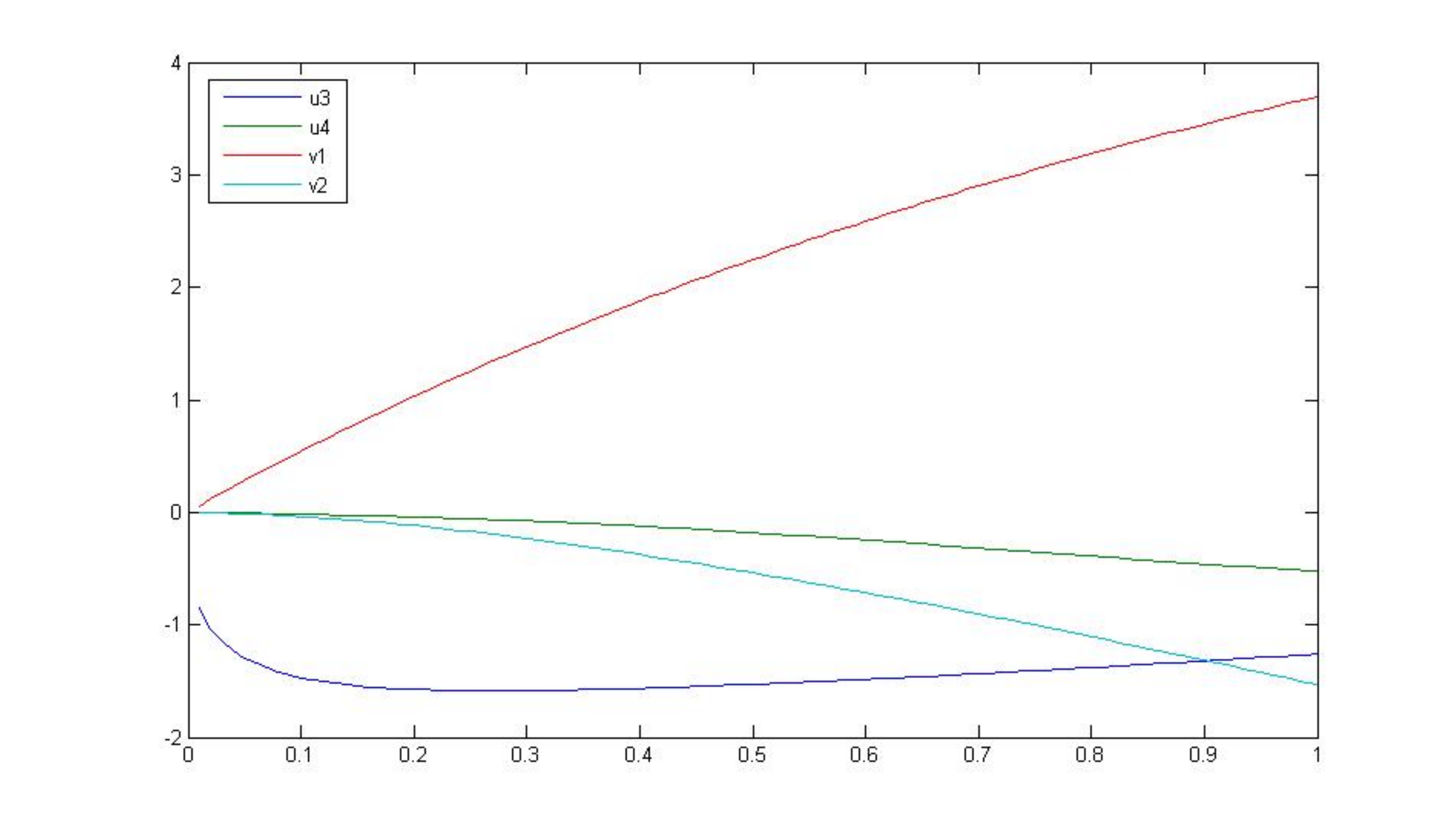}}\caption{The estimated values of $u_3(0;m)$, $u_4(0;m)$, $v_1(0;m)$, and $v_2(0;m)$ for the nonzero components of $\gamma(0;m)$ for $0<m\leq 1$.}\label{figure3}
\end{figure}

\begin{figure}\scalebox{0.4}{\includegraphics{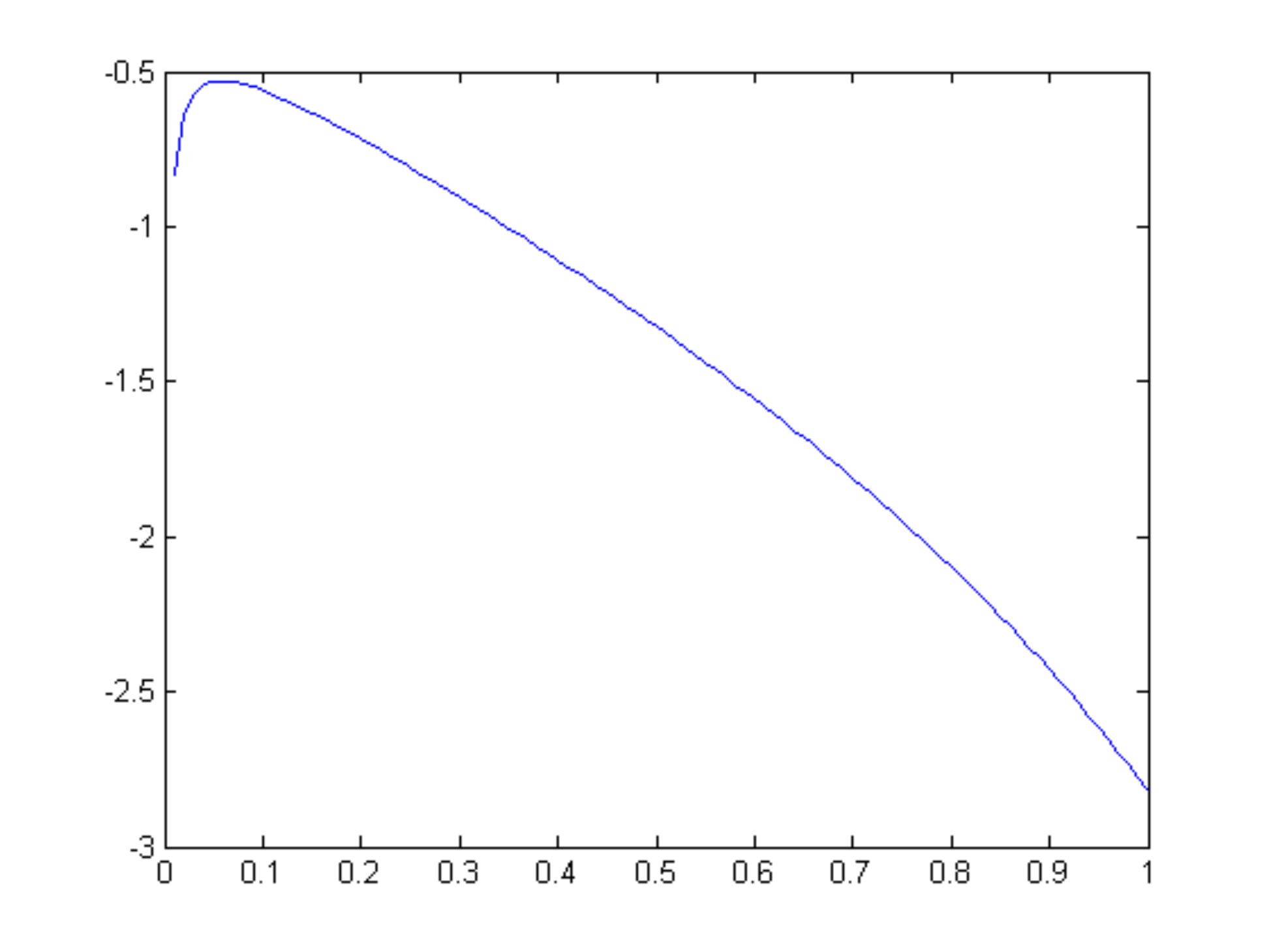}}\caption{The energy $\hat E(m)$ of the periodic orbits $\gamma(s;m)$ for $0<m\leq 1$.}\label{figure4}
\end{figure}

\begin{figure}\scalebox{0.4}{\includegraphics{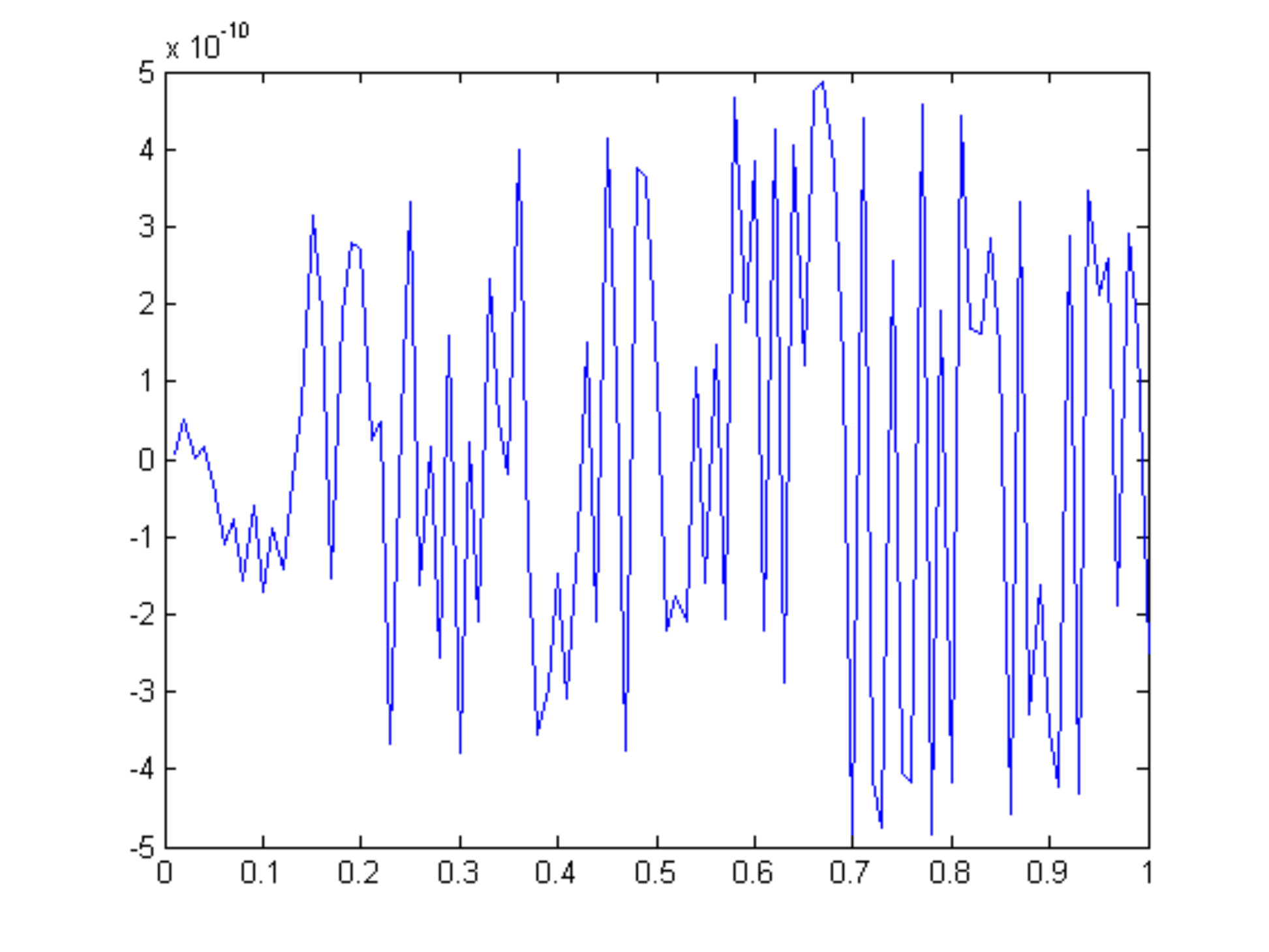}}\caption{The value of $\hat\Gamma(\gamma(\pi/4;m))$ for the optimized trigonometric polynomial approximation of $\gamma(s;m)$ over $0<m\leq 1$.}\label{figure5}
\end{figure}

This numerical method only works well if we have a good initial guesses for $a_i$, $b_i$, $c_i$, $d_i$, $i=1,\dots,5$, and $\hat{E}(m)$ for some value of $m$. This we have when $m=1$. We use our estimate of $\gamma(s;1)$ and $\hat E(1)$ to provide the initial guesses for $a_i$, $b_i$, $c_i$, $d_i$, $i=1,\dots,5$, and $\hat{E}(m)$ for $m=0.99$. The numerical algorithm produces a trigonometric polynomial approximation of $\gamma(s;0.99)$ and an estimate of $\hat E(0.99)$, which then provide the initial guesses for $a_i$, $b_i$, $c_i$, $d_i$, $i=1,\dots,5$, and $\hat{E}(m)$ for $m=0.98$. We continue decreasing $m$ by $0.01$ and using the numerical algorithm until we reach $m=0.01$. In Figure \ref{figure3}, we plot the graphs of the numerical estimates of $u_3(0;m)$, $u_4(0;m)$, $v_1(0;m)$, and $v_2(0;m)$. Notice that $v_1(0;m)$ and $v_2(0;m)$ satisfy Equation (\ref{rmomentum}), as is expected from the regularization of the simultaneous binary collisions. In Figure \ref{figure4}, we plot the graph of the numerical estimate of $\hat E(m)$. In Figure \ref{figure5}, we graph of value of $\hat\Gamma(\gamma(\pi/4;m))$ over $0<m\leq 1$.

We use the optimized trigonometric approximations of $\gamma(s;m)$, $0<m<1$, to estimate $\gamma(-\pi/4;m)$, which are the conditions corresponding to $s=0$ before the shift of $s$.  For each $0<m<1$, the conditions $\gamma(-\pi/4;m)$ satisfy Lemma \ref{initial}, and so they corresponds to the initial conditions of interest as given in the Introduction.

\begin{figure}\scalebox{0.5}{\includegraphics{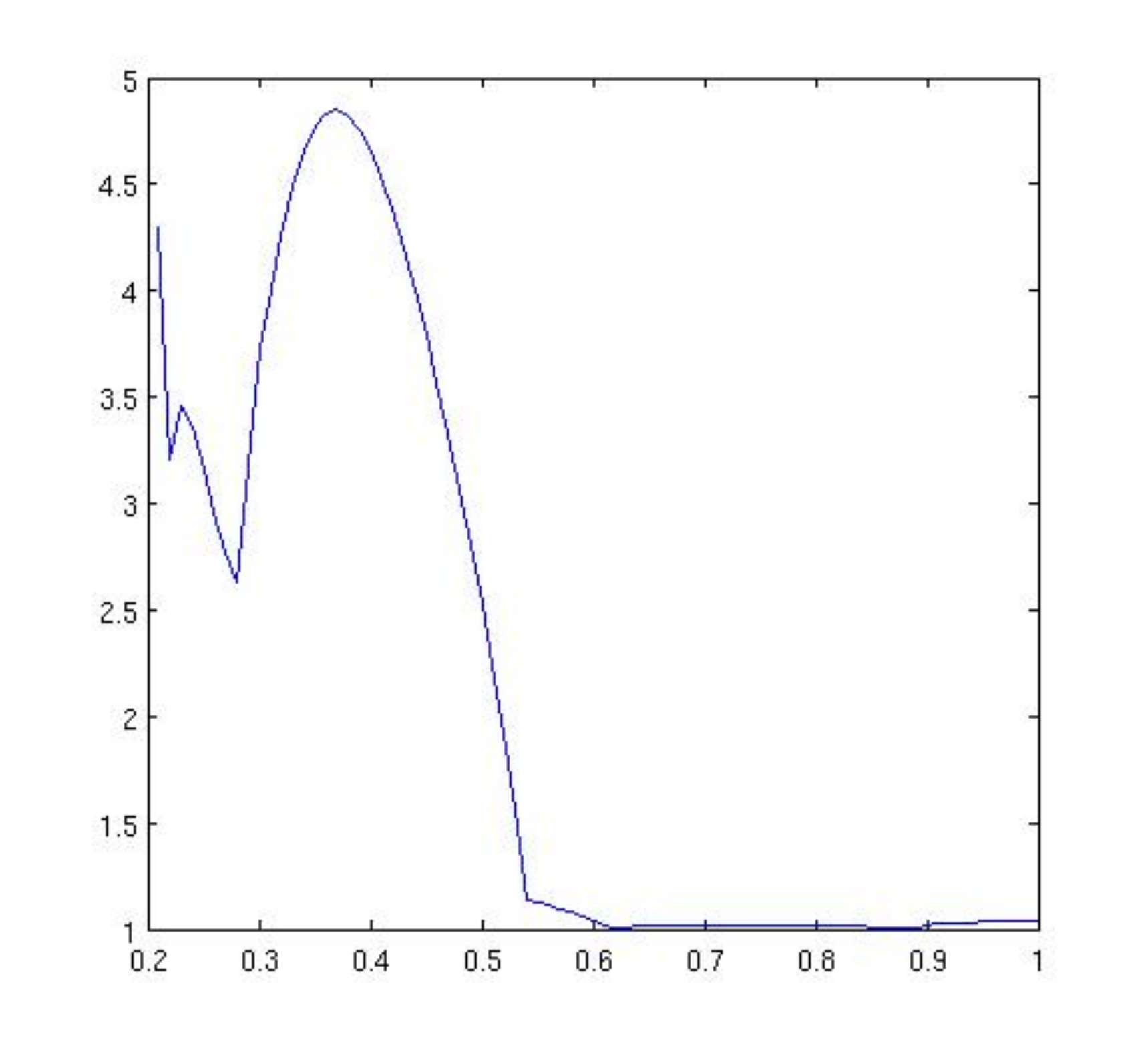}}\caption{The maximum modulus of the eight eigenvalues of $X(s;m)$ for $0.2<m\leq 1$.}\label{figure6}\end{figure}

We numerically estimate the monodromy matrices $X(s;m)$ and their eigenvalues for $\gamma(s;m)$ with $0<m\leq 1$. We used a Runge-Kutta order $4$ algorithm coded in MATLAB with a fixed time step of $2\pi/50000$ at increments of $0.01$ in $m$. In Figure \ref{figure6} we graph the maximum modulus of the eight eigenvalues of $X(s;m)$ over $0.2<m\leq1$. The maximum modulus of the eigenvalues of $X(s;m)$ for $0<m\leq 0.2$ becomes large as $m\to 0$, and are therefore not included in Figure \ref{figure6}. The maximum modulus for the eight eigenvalues shows that $\gamma(s;m)$ is linearly unstable for $0<m\leq0.53$. The error in the numerical estimation of the repeated eigenvalue $1$ of $X(s;m)$ is large enough to account for the graph in Figure \ref{figure6} not lying at $1$ over $0.54\leq m\leq 1$. However, for $0.54\leq m\leq 1$, the other six eigenvalues come in three distinct complex conjugate pairs that are each of modulus one and located away from $\pm 1$. This shows that $\gamma(s;m)$ is linearly stable for $0.54\leq m\leq 1$.

\vskip0.2cm

\noindent{\bf Acknowledgements}. We thank the referees for their valuable comments that improved the paper.

\end{document}